\documentclass{article}
\usepackage{array}
\usepackage[utf8]{inputenc}
\usepackage{microtype}
\usepackage{graphicx}
\usepackage{subcaption}
\usepackage{caption}
\usepackage{booktabs} 
\usepackage[margin=1.0in]{geometry}
\usepackage{mathtools}
\usepackage{amsfonts}           
\usepackage{amsthm}
\usepackage{amsmath}
\usepackage{amssymb}
\usepackage{stmaryrd}
\usepackage{xcolor}
\usepackage{caption}
\usepackage{algorithm}
\usepackage{algpseudocode}
\usepackage{comment}
\usepackage{authblk}
\usepackage{hyperref}
\usepackage{multirow}
\usepackage[table]{xcolor}
\usepackage{hyperref}

\newcommand{\email}[1]{\href{mailto:#1}{\textit{#1}}}

\begin{document}

\title{A Persistent Homology Signature of Knotting
}

\author[1,2,*]{Aurelie Jodelle Kemme}
\author[1]{Collins A.\ Agyingi}
\author[4]{Colleen Farrelly}
\author[3,*]{Agnese Barbensi}
\affil[1]{\small Department of Mathematical Sciences, University of South Africa, Pretoria, South Africa}
\affil[2]{\small African Institute for Mathematical Sciences, Research and Innovation Centre, Kigali, Rwanda}
\affil[3]{\small School of Mathematics and Physics, the University of Queensland, Brisbane, Australia}
\affil[4]{\small Post Urban Ventures, London, UK}
\affil[*]{Corresponding author: ~\email{a.barbensi@uq.edu.au}}

\maketitle

\abstract{}
We ask whether knotting can be recognised using persistent homology. Starting from a point-cloud representation of a curve, we compute one-dimensional persistent homology, extract cycle representatives, and assign a hypergraph curvature-based score to these cycles. Motivated by proteins but tested more broadly, the method reveals systematic differences between knotted and unknotted structures in both protein families and synthetic examples. This suggests that knotting leaves a detectable persistent-homology-based signature.
\section{Introduction}\label{sec1}

\begin{figure}[ht]
  \centering
  \includegraphics[width=0.5\textwidth]{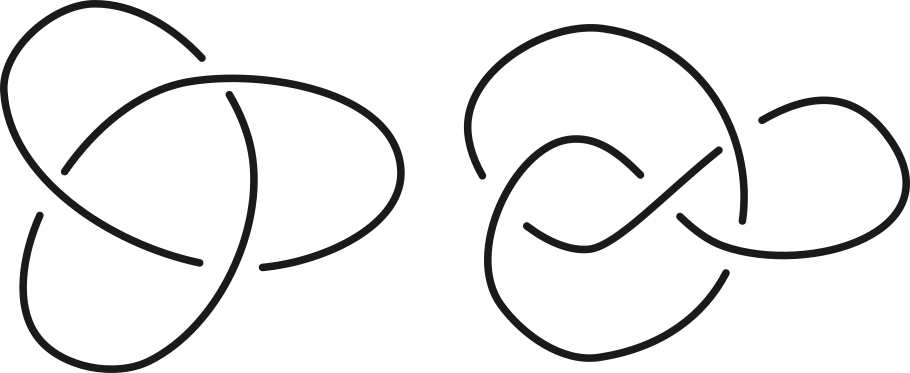}
  \caption{On the left-hand side, the trefoil
    knot ($3_1$), the simplest non-trivial knot, and the most commonly observed knot in proteins. On the right-hand side, the figure-eight knot ($4_1$), the second simplest knot.}
  \label{fig:knots}
\end{figure}

Knotted proteins are proteins whose backbone chain traces a non-trivial topological knot in three-dimensional space. They were first identified systematically by Mansfield~\cite{mansfield1994there}, and later confirmed across a growing number of protein families~\cite{taylor2000deeply,king2007identification,faisca2015knotted}. The KnotProt 2.0 database~\cite{jamroz2015knotprot, dabrowski2019knotprot} now catalogues these structures systematically, spanning knot types from the trefoil ($3_1$) to six-crossing topologies. The most common knotted topologies are the trefoil ($3_1$) knot (Fig.~\ref{fig:knots}, left-hand side) and the figure-eight ($4_1$) knot (Fig.~\ref{fig:knots}, right-hand side).

The presence of knots in proteins challenges standard models of protein folding. Despite decades of research, a full understanding of \textit{why} and \textit{how} proteins fold to create non-trivial knots is still out of reach~\cite{tubiana2024topology}. An important question in knotted-protein studies is whether there are recurring structural motifs that distinguish knotted proteins from unknotted ones. One possible approach is to look for \emph{knot-promoting} structures, \textit{i.e.}~local geometric features that enable (or stabilise) the strand passage needed to create an entanglement. For example, earlier work identified \emph{knot-promoting loops} in several families: short loop insertions whose removal (or virtual bridging) eliminates the knot~\cite{potestio2010knotted}.

An interesting and well-studied example is provided by trefoil-knotted AOTCases and their unknotted OTCase counterparts. Proteins in these families have almost superimposable structures, the key difference being a localised strand passage deep in the fold, which creates the knot in the AOTCase case~\cite{virnau2006intricate,potestio2010knotted}. This makes the pair a natural sandbox to ask whether knotting leaves a detectable, local topological signature beyond standard geometric descriptors.

Topological data analysis (TDA) provides a natural framework to search for this type of structural signal. Persistent homology (PH), one of the main tools in TDA, is particularly suited to detect non-trivial structural features in complex geometric data~\cite{ghrist2008barcodes}. In its simplest formulation, PH starts from a set of points in Euclidean space and tracks how topological features appear and disappear as a scale parameter grows. The output can be represented as barcodes or persistence diagrams, giving a multiscale topological summary of the data~\cite{ghrist2008barcodes}. Persistent-homology methods have already been extensively used in protein studies, \textit{e.g.}~to describe folding dynamics and conformational changes, and to extract informative structural fingerprints~\cite{xia2014persistent,ichinomiya2020protein,kovacevnikolic2016binding,bramer2020atom}. More recently, PH pipelines that also exploit explicit cycle representatives have made the output more local and interpretable, and have enabled analysis at scale across large protein datasets~\cite{barbensi2022hypergraphs,madsen2025proteinuniverse}.

In the AOTCases \textit{vs} OTCases context, persistent homology has been used to detect a local structural signal present in knotted AOTCases and absent in unknotted OTCases~\cite{benjamin2023homology}. Specifically, the knotted structures exhibit an additional topological feature: an extra $H_1$ class appears in their persistent homology compared to the unknotted homologues, and cycle representatives localise this difference to the region responsible for the strand passage. This provides a direct and interpretable link between persistent homology and protein knotting. Related work also shows that PH captures geometric information in knotted embeddings more generally, for example in the study of knot confinement~\cite{celoria2022knotconfinement}. \\

Motivated by these results, we ask whether the same idea extends more generally. In particular, we investigate whether persistent homology can recognise knottedness in proteins through simple and robust structural signatures, beyond the AOTCase/OTCase case study. More broadly, we ask whether knotting can be detected from PH-based descriptors across different families of knotted proteins and knotted curves. 

To test this idea, we build on the persistent-homology pipeline introduced in Benjamin~\textit{et al.}~\cite{benjamin2023homology} and developed further in the hyperTDA framework~\cite{barbensi2022hypergraphs}. We compute one-dimensional PH of the protein backbone point cloud and extract explicit $H_1$ cycle representatives. As in hyperTDA, we then encode these representatives as hyperedges of a PH-hypergraph~\cite{barbensi2022hypergraphs}. We assign to each hyperedge the unweighted undirected Forman-Ricci curvature~\cite{leal2021forman}, obtaining a simple local scalar descriptor of how strongly the residues in one cycle are shared with other cycles. \\

\begin{table}[ht]
\centering
\caption{Representative examples of knotted and slipknotted proteins together with their unknotted homologs sharing more than 40\% sequence similarity, as reported in KnotProt. Here ``Length'' corresponds to the sequence length, defined as the number of C$_\alpha$ atoms in the protein structure. Unknotted homologs are proteins with similar sequences but without knotting. The ``PDB ID'' is the identifier assigned by the Protein Data Bank.}
\label{tab:dataset_size}
\begin{tabular}{l l >{\centering\arraybackslash}p{1.6cm} c r  l}
\toprule
\textbf{Protein Family} & \textbf{PDB ID} & \textbf{Slip/Knot} & \textbf{Knot Type} & \textbf{Length} & \textbf{Class} \\
\midrule
 Permease & 6kkt & \cellcolor{purple!50} Slipknot & S4(1) & 537 & 7tti \\
Solute symporter family & 2xq2 & \cellcolor{purple!50} Slipknot & S4(1) & 538 & 7v19 \\
\cmidrule(lr){2-6}
Myxobacterial phytochrome & 6bay & \cellcolor{gray!20} Knot & K4(1) & 505 &  4r70 \\
Reductoisomerase           & 5e4r & \cellcolor{gray!20} Knot & K4(1) & 466 &  5yeq \\
Bacteriophytochrome        & 5i5l & \cellcolor{gray!20} Knot & K4(1) & 480 &  4s21 \\
\cmidrule(lr){2-6}
\multirow{3}{*}{Carbonic anhydrase} 
 & 6rqn & \cellcolor{gray!20} Knot & K+3(1) & 256 & 6rqw \\
 & 3f7b & \cellcolor{gray!20} Knot & K+3(1) & 257 &5jn9 \\
 & 6r71 & \cellcolor{purple!50} Slipknot & S+3(1) & 261  & 4ww8 \\
\cmidrule(lr){2-6}
Carbamoyltransferase       & 2yfk & \cellcolor{gray!20} Knot & K+3(1) & 385 & 6jwx \\
\cmidrule(lr){2-6}
Alkaline phosphatase        & 3wbh & \cellcolor{purple!50} Slipknot & S+3(1) & 497  & 6qsq \\
Phosphoglycerate mutase     & 5kgn & \cellcolor{purple!50} Slipknot & S+3(1) & 530 & 7kng \\
\bottomrule
\end{tabular}
\end{table}

Proteins knot topologies are determined computationally by artificially \textit{closing} the space curve to form a closed loop. These closures, often performed stochastically,  characterise the space curve by the knot type it assumes upon closure. In this context, a protein is called  \textit{slipknotted} if it admist a knotted subchain, but the overall backbone is classified as unknotted. We apply this curvature-based pipeline to knotted proteins and their homologous unknotted counterparts across four familisies: knotted and slipknotted positive trefoils and figure eight knots K+3(1), S+3(1), K4(1), S4(1), see Table~\ref{tab:dataset_size}. We further test the trend on randomly generated knotted curves with increasing length.\\

Across all families examined, we observe a consistent shift of the hyperedge-curvature distribution towards more negative values in knotted structures compared to unknotted ones, suggesting that strong negative undirected Forman-Ricci curvature captures a robust signature of chain entanglement, coherently with the intuitive idea that geoemetric sub-structures (and thus persistent cycles) are more interwined in the presence of knots. \\

This paper asks whether PH can recognise protein knottedness, and shows that a very simple PH-based curvature descriptor does so consistently across several families. So while the computation is simple, the main point is new: knotting appears to leave a detectable and quantifiable signature in PH. This perspective is separate from existing connections between persistent homology and curvature in other settings~\cite{bubenik2020persistent}. Here, the connection is realised through cycle-derived hyperedges rather than curvature-based filtrations on graphs. \\

The paper is organised as follows. Section~2 describes the protein datasets, the computational pipeline, and the statistical comparison strategy. Section~3 presents the results for each of the four protein families, followed by a validation study on synthetic polymer loops. Section~4 provides a conclusion and discussion of biological implications, limitations, and directions for future work.

\section{Methods}

This section contains the description of the datasets considered and of the mathematical framework and computational implementation underlying our analysis.

\subsection{Datasets}
\begin{enumerate}
    \item \textbf{Protein Dataset} The protein dataset is restricted to positive trefoil knots $K{+}3(1)$, figure-eight knots $K4(1)$, and their respective slipknot counterparts $S{+}3(1)$ and $S4(1)$, as they are the most extensively catalogued knot types in the database. Several additional classes, $K{-}3(1)$ (sample size $n=33$), $S{-}3(1)$ ($n=38$), $K{-}5(2)$ ($n=26$), and $K6(1)$ ($n=3$) were excluded on account of insufficient sample sizes to support reliable distributional comparisons. 

    For each knotted class, we extract homologous proteins from the \textit{Similar chains (by sequence)} entries in KnotProt, which correspond to a sequence identity threshold of $\geq 40\%$~\cite{dabrowski2019knotprot}. We then partition all proteins considered into \textit{homology classes}, and homologs in each class into either \textit{knotted} or \textit{unknotted}, depending on if they form knots, as detailed in Table~\ref{tab:dataset_size}, see also the Data accessibility Section. Table~\ref{tab:dataset_methods_compact2} reports the final sample sizes for each class. For each protein structure, we consider \textit{xyz}-coordinates of the $\alpha$-carbon ($C_\alpha$) building the protein backbone, yelding one point cloud per protein.  These are exctracted from Knotprot~\cite{dabrowski2019knotprot} or the Protein Data Bank (PDB)~\cite{berman2000protein}. 

    \item []

    \item \textbf{Random Dataset} While this study is inspired by knotted proteins, our aims was to show weather persistent homology could detect knottinnes in full generality. To test this idea, we analyse random knotted loops of varying lengths. These are generated using Pyhton's \texttt{Topoly} package~\cite{dabrowski2021topoly}. In \texttt{Topoly}, loops are generated as self-avoiding random walks. This procedure yields piece-wise linear closed curves of a prescribed length. We classify their topology by computing the Jones polynomial, and labelling a loop as unknotted when this is $1$, and as knotted otherwise. Generation was performed over nine chain lengths from $100$ to $500$ with steps of $50$. For each length, we generate $500$ knotted loops and $500$ unknotted ones by rejection sampling. 
\end{enumerate}

\begin{table}
\centering
\footnotesize
\caption{Summary of dataset sizes and chain-length distributions. The ``Size'' column gives the number of unique protein structures (or loops) in each class. ``Median Chain length'' is the median number of $C_\alpha$ atoms across all proteins. For polymer loops, chain length is fixed by design at $L \in \{100, \ldots, 500\}$.}
\label{tab:dataset_methods_compact2}
\begin{tabular}{l l r l r}
\toprule
\textbf{(Slip)Knot Type} & \textbf{Slip/Knot} & \textbf{Size} & \textbf{Source} & \textbf{Median chain length} \\
\midrule
$S4(1)$  & \cellcolor{purple!20}Slipknot & $127$ & KnotProt & $507$  \\
Unknotted homologs $S4(1)$ & \cellcolor{gray!20}unknotted    & $48$  & PDB & $537$  \\
\midrule
$K4(1)$  & \cellcolor{purple!20}Knot     & $70$  & KnotProt & $478$  \\
Unknotted homologs $K4(1)$ & \cellcolor{gray!20}unknotted  & $13$  & PDB & $462$  \\
\midrule
$K{+}3(1)$  & \cellcolor{purple!20}Knot     & $615$ & KnotProt & $257$\\
Unknotted homologs $K{+}3(1)$ & \cellcolor{gray!20}unknotted  & $94$  & PDB & $310$\\
\midrule
$S{+}3(1)$ & \cellcolor{purple!20}Slipknot & 1{,}211 & KnotProt & $259$\\
Unknotted homologs $S{+}3(1)$ & \cellcolor{gray!20}unknotted  & $340$  & PDB & $427$\\
\midrule
Knotted polymer loops & \cellcolor{purple!20}mixed   & 500 & \texttt{Topoly} & $100$-$500$ \\
Unknotted polymer loops & \cellcolor{gray!20}unknot    & 500 & \texttt{Topoly} & $100$-$500$ \\
\bottomrule
\end{tabular}
\end{table}

\subsection{The computational pipeline}
\begin{figure}[ht]
  \centering
  \includegraphics[width=0.95\linewidth]{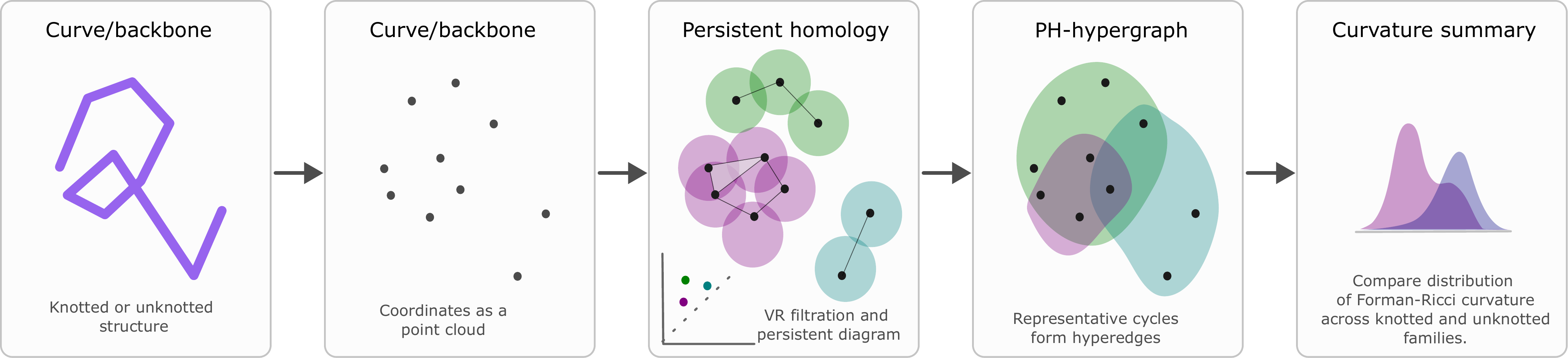}
  \caption{Overview of the computational pipeline. From left to right: from a spatial curve, to a point cloud. One-dimensional persistent homology is computed via a Vietoris-Rips filtration; cycle representatives of each persistent homology class are extracted using the involuted algorithm of \texttt{Ripserer.jl}; these representatives are used as hyperedges to form the PH-hypergraph (note that the hypergraph shown here is a schematic sketch and not a true representation of the cycle representatives of the actual data); and the unweighted undirected Forman-Ricci curvature is computed on each hyperedge. The distribution of median curvature is compared across knotted and unknotted families.}
  \label{fig:pipeline}
\end{figure}

The pipeline we propose (Fig.~\ref{fig:pipeline}) builds on the PH-hypergraph construction of~\cite{barbensi2022hypergraphs}. Given a protein chain with $N$ residues, we consider the set of $C_\alpha$ coordinates $\mathcal{X}\subset \mathbb{R}^3$, indexed by residue number. Likewise, for a randomly generated piecewise-linear loop of length $N$, the point cloud $\mathcal{X}\subset \mathbb{R}^3$ is given by the coordinates of its $N$ vertices. 

\subsubsection{Persistent Homology via Vietoris-Rips Filtration}
For each point cloud $\mathcal{X}$, we compute persistent homology using the Vietoris--Rips filtration over $\mathbb{Z}_2$, implemented in \texttt{Ripserer.jl}~\cite{cufar2020ripserer,bauer2021ripser}. In this setting, persistent homology tracks how loops appear and disappear as the scale parameter increases. More precisely, at each scale $\varepsilon$, the Vietoris-Rips complex connects points whose pairwise distance is at most $\varepsilon$, and fills in higher-dimensional simplices whenever all corresponding lower-dimensional faces are present. One-dimensional persistent classes therefore detect loop-like structures in the point cloud that persist across a range of scales. Since our interest is precisely in these loop-like features, we restrict attention to degree one. This computation yields two primary outputs for each protein/loop: (i)~a set of cycle representatives, each representing a one-dimensional persistent classes, encoded as a collection of one-simplices (edges) where each index refers to a residue position in the backbone; and (ii)~a persistence diagram recording the birth and death scales of each feature. Long-lived cycles are expected to correspond to structurally significant regions.\\

For each persistent class $c$, let $\mathrm{rep}(c)$ denote a chosen cycle representative, written as a collection of edges. We then define the corresponding hyperedge by
\begin{equation}
    \sigma_c=\bigcup_{[v_1,v_2]\in \mathrm{rep}(c)} \{v_1,v_2\}\subseteq \{1,\dots,N\}.
\end{equation}
This yields a PH-hypergraph $\mathcal{H}=(\mathcal{V},\mathcal{E})$, where $\mathcal{V}=\{1,\dots,N\}$ and $\mathcal{E}=\{\sigma_c\}$. Its incidence matrix is
\begin{equation}
    \mathcal{B}_{i,\sigma_c}=
    \begin{cases}
        1, & \text{if } i\in \sigma_c,\\
        0, & \text{otherwise.}
    \end{cases}
\end{equation}
Thus each hyperedge is the support of one persistent cycle representative, and overlaps between hyperedges record vertices shared by different cycles.

A standard point to note is that the persistent classes are well defined, but their cycle representatives are far from being unique. Different matrix reductions, tie-breaking rules, or optimisation criteria may produce different generators for the same interval in the barcode. 

This issue is well known whenever representatives are used geometrically rather than only through their birth-death pairs~\cite{zhang2025topological, li2021minimal}. In our setting, the PH-hypergraph is therefore defined with respect to a chosen family of representatives, namely those returned by \texttt{Ripserer.jl}. This gives a consistent construction across the dataset, while leaving open the question of how alternative choices of generators, for instance shorter or more localised representatives, might affect the resulting hyperedges geometry and curvature values. However, as detailed in~\cite{barbensi2022hypergraphs}, hypergraph analysis is experimentally shown to be robust by different choices of cycles.

\subsubsection{Forman-Ricci Curvature on the PH-Hypergraph}

The first departure from the PH-hypergraph analysis in Barbensi \textit{et al}~\cite{barbensi2022hypergraphs} is our choice of descriptor. Rather than applying centrality measures or community detection, we compute the \emph{unweighted undirected Forman-Ricci curvature} on each hyperedge, following Leal~\textit{et al.}~\cite{leal2021forman}. For a hyperedge $e \in \mathcal{E}$, the general Forman-Ricci curvature is:
\begin{equation}\label{ricci-curvature}
\mathrm{F}(e) = w_{e}\left[\sum_{k \in e} \frac{w_k}{w_e} - \sum_{\substack{e_i \ni k \\ e_i \neq e}} \frac{w_k}{\sqrt{w_e\, w_{e_i}}}\right],
\end{equation}
where $w_k$ and $w_e$ denote the weights of vertex $k$ and hyperedge $e$ respectively, and the second sum runs over all hyperedges $e_i$ incident to any vertex $k \in e$, excluding $e$ itself. In the unweighted setting ($w_k = w_e = 1$), this simplifies to:
\begin{equation}
    \mathrm{F}(e) = 2|e| - \sum_{k \in e} d_k \;=\; 2|e| - D,
\end{equation}
where $|e|$ is the cardinality of hyperedge $e$, $d_k$ is the degree of vertex $k$ (the number of hyperedges in $\mathcal{E}$ containing $k$), and $D = \sum_{k \in e} d_k$ is the total degree sum over the hyperedge. For fixed $|e|$, the curvature decreases as the constituent vertices participate in more hyperedges, directly encoding the degree of topological entanglement at each cycle. This interpretable local geometric quantity provides a descriptor at the level of individual hyperedges, rather than a network-level ranking or partition, which is the key motivation for preferring it over centrality-based alternatives. To each protein, we associate a distribution of hyperedge curvature values, one per $H_1$ cycle. The median is used as the scalar summary because it is resistant to extreme curvature values at either tail; a small number of isolated hyperedges with very high or very low $F(e)$ would disproportionately influence the mean. The median captures the \emph{global} curvature regime of the hypergraph, which is the quantity of biological interest. In this study, per-protein median is used in all subsequent statistical comparisons.

\section{Results}

In this framework, each $H_1$ persistent cycle is encoded as a hyperedge with Forman-Ricci curvature $F(e) = 2|e| - D$, where $|e|$ is the number of residues/vertices in the hyperedge and $D = \sum_{k \in e} d_k$ is the total degree of those vertices across all hyperedges. \\

We propose the following hypothesis: \textit{knotting manifests as an increased overlap between persistent cycles, raising their individual degrees and therefore increasing $D$. Since $D$ enters the curvature formula with a negative sign, any such increase in $D$ directly predicts more negative $F(e)$}. \\

Knotted curves are therefore hypothesised to exhibit, on average, more negative median Forman-Ricci curvature than their unknotted counterparts, with reduced variance, because the entanglement constraint narrows the achievable range of curvature values. This hypothesis is directly testable and is not verified at the level of individual residue localisation within this study; we test it at the distributional level across four protein families and validate it on synthetic loops in Section~\ref{subsec:loops}.

\subsection{Trefoil families: \texorpdfstring{$K{+}3(1)$}{K+31} and \texorpdfstring{$S{+}3(1)$}{S+31} versus unknotted homologs}
\label{subsec:trefoil_families}

We begin with the trefoil families, comparing the knotted family $K{+}3(1)$ and the slipknotted family $S{+}3(1)$ against their respective unknotted homologs. In both cases, the same overall trend is observed: knotted or slipknotted proteins have curvature distributions that are shifted towards more negative values and are more concentrated than those of their unknotted comparators.\\

\begin{figure}[ht]
    \centering
    \begin{subfigure}[b]{0.45\textwidth}
        \centering
        \includegraphics[width=\linewidth]{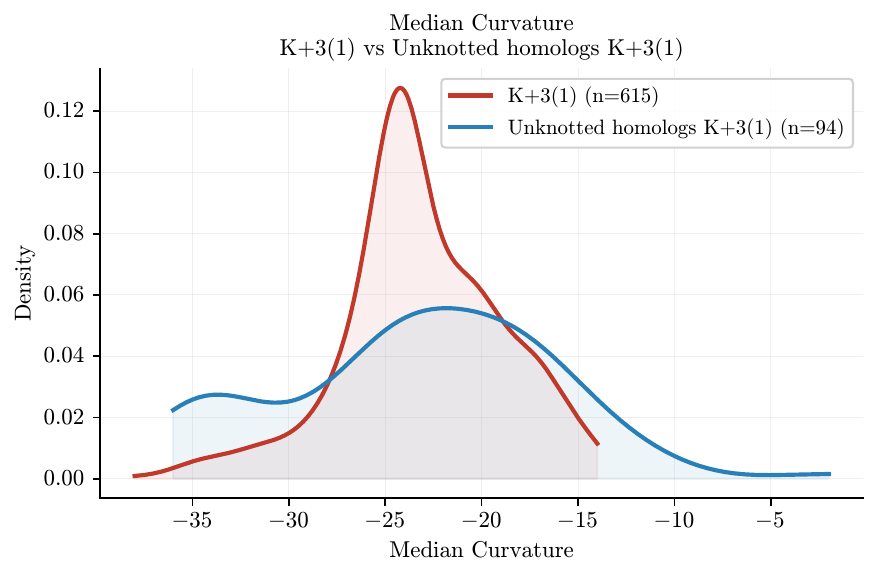}
        \caption{}
        \label{fig:K+3(1)_Curv_kde_median}
    \end{subfigure}
    \hfill
    \begin{subfigure}[b]{0.45\textwidth}
        \centering
        \includegraphics[width=\linewidth]{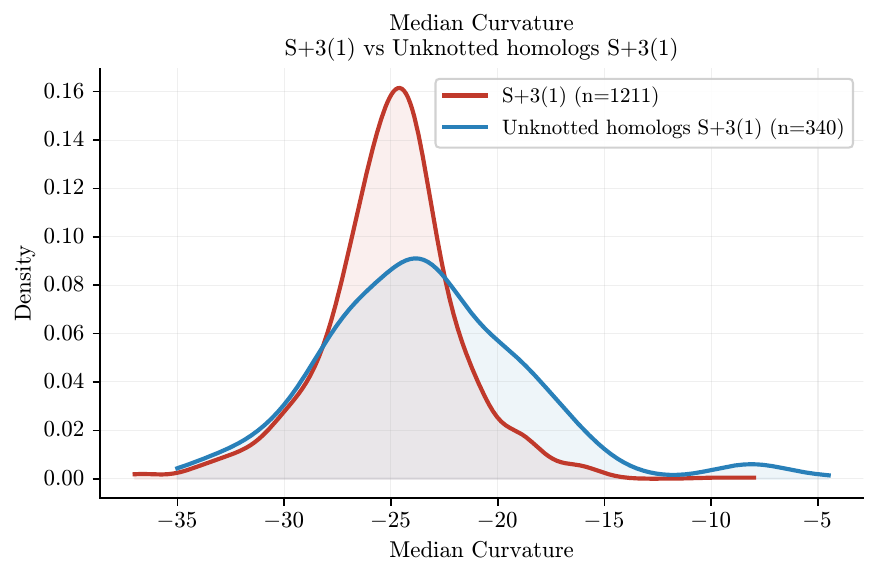}
        \caption{}
        \label{fig:S+3(1)_Curv_kde_median}
    \end{subfigure}
    \caption{Kernel Density Estimation (KDE) of the per-protein median undirected Forman--Ricci curvature for the trefoil families and their unknotted homologs.
    \textbf{(a)}~Trefoil knot proteins ($K{+}3(1)$, $n=615$) versus unknotted homologs ($n=94$). The knotted family is shifted towards more negative curvature values and is less dispersed.
    \textbf{(b)}~Trefoil slipknot proteins ($S{+}3(1)$, $n=1{,}211$) versus unknotted homologs ($n=340$). Both groups peak near $-25$, but the slipknotted family is more concentrated, while the unknotted homologs are broader and extend further towards less negative curvature.}
    \label{fig:trefoil_families_kde}
\end{figure}

For the right-handed trefoil knot $K{+}3(1)$ ($n=615$), compared against $94$ unique unknotted homologs, the Kernel Density Estimation (KDE) in Fig.~\ref{fig:trefoil_families_kde}\subref{fig:K+3(1)_Curv_kde_median} shows a clear left shift. The $K{+}3(1)$ distribution peaks near $-25$ and is more compact, whereas the unknotted homologs peak near $-19$ and are more dispersed. The group medians are $-24$ ($K{+}3(1)$) and $-23$ (unknotted). The standard deviations (SDs) are $4.13$ and $7.00$, respectively, giving a variance ratio of $0.350$. This difference is supported statistically: the Kolmogorov-Smirnov (KS) test gives $D = 0.167$ with false discovery rate (FDR)-corrected $p_{\mathrm{FDR}} = 0.025$ (Table~\ref{tab:ks_test}), and the Levene test for equality of variances gives $F = 51.19$ with $p_{\mathrm{FDR}} = 4.2 \times 10^{-12}$ (Table~\ref{tab:levene_test}).
The same qualitative picture appears for the trefoil slipknot $S{+}3(1)$, the largest family in the dataset ($n=1{,}211$), compared against $340$ unknotted homologs. As shown in Fig.~\ref{fig:trefoil_families_kde}\subref{fig:S+3(1)_Curv_kde_median}, both groups peak near $-25$, but the $S{+}3(1)$ family is much more concentrated, while the unknotted homologs are broader and extend further towards less negative values, with a secondary peak near $-14$. The group medians are $-25$ ($S{+}3(1)$) and $-24$ (unknotted). The corresponding SDs are $3.18$ and $4.85$, giving a variance ratio of $0.432$. Both the variance difference and the overall distributional difference are significant after FDR correction: the Levene test gives $F = 75.53$ with $p_{\mathrm{FDR}} = 3.6 \times 10^{-17}$, and the KS test gives $D = 0.214$ with $p_{\mathrm{FDR}} = 1.6 \times 10^{-10}$ (Tables~\ref{tab:levene_test} and~\ref{tab:ks_test}).

Taken together, these two trefoil families show a consistent tendency for knotted or slipknotted proteins to occupy a more negative and less variable curvature regime than their unknotted homologs.

\subsection{Figure-eight families: \texorpdfstring{$K4(1)$}{K4(1)} and \texorpdfstring{$S4(1)$}{S4(1)} versus unknotted homologs}
\label{subsec:figure_eight_families}

We next consider the figure-eight families. As in the trefoil case, both the knotted family $K4(1)$ and the slipknotted family $S4(1)$ show a tendency towards more negative curvature values than their unknotted homologs, although the statistical strength of the effect differs between the two comparisons.

\begin{figure}[ht]
    \centering
    \begin{subfigure}[b]{0.45\textwidth}
        \centering
        \includegraphics[width=\linewidth]{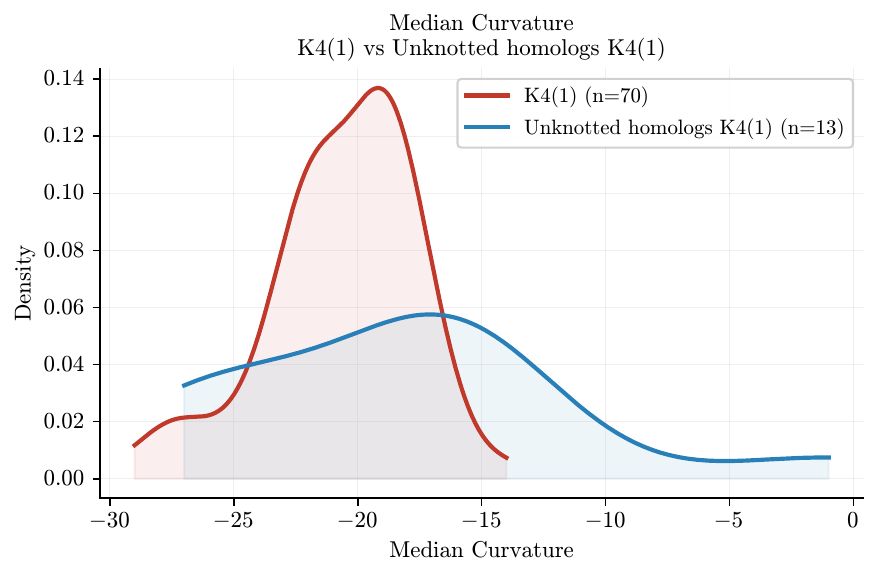}
        \caption{}
        \label{fig:K4(1)_Curv_kde_median}
    \end{subfigure}
    \hfill
    \begin{subfigure}[b]{0.45\textwidth}
        \centering
        \includegraphics[width=\linewidth]{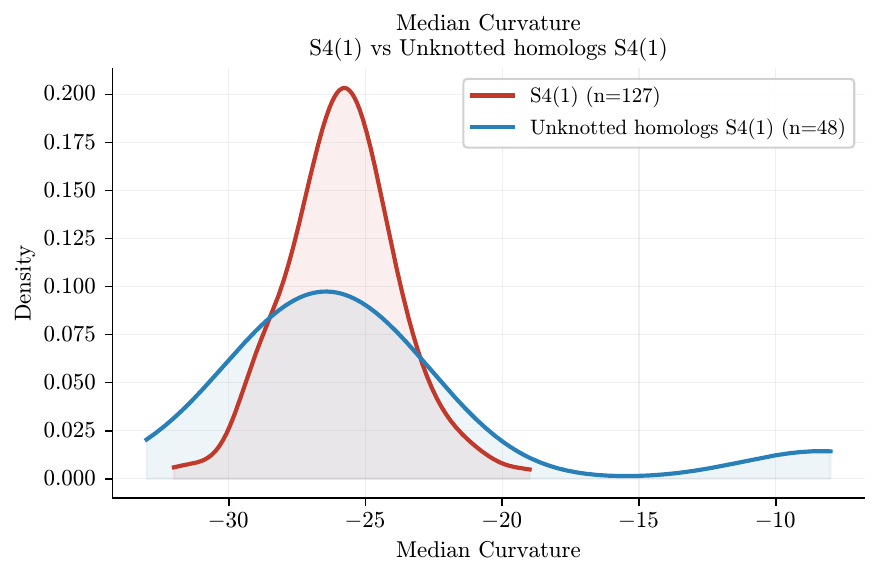}
        \caption{}
        \label{fig:S4(1)_Curv_kde_median}
    \end{subfigure}
    \caption{Kernel Density Estimation (KDE) of the per-protein median undirected Forman--Ricci curvature for the figure-eight families and their unknotted homologs.
    \textbf{(a)}~Figure-eight knot proteins ($K4(1)$, $n=70$) versus unknotted homologs ($n=13$). The knotted family is shifted towards more negative curvature values and is substantially less dispersed.
    \textbf{(b)}~Figure-eight slipknot proteins ($S4(1)$, $n=127$) versus unknotted homologs ($n=48$). The slipknotted family is narrower and more concentrated, although the group medians coincide.}
    \label{fig:figure_eight_families_kde}
\end{figure}

For the figure-eight knot $K4(1)$ ($n=70$), compared against $13$ unknotted homologs, the KDE in Fig.~\ref{fig:figure_eight_families_kde}\subref{fig:K4(1)_Curv_kde_median} shows that the $K4(1)$ family is shifted towards more negative curvature values. The $K4(1)$ distribution is concentrated between $-25$ and $-15$ with a sharp peak near $-20$, whereas the unknotted homologs are broader, right-skewed, and peak near $-17$. The group medians are $-20$ ($K4(1)$) and $-17$ (unknotted). The corresponding SDs are $2.93$ and $6.96$, giving a variance ratio of $0.189$. Both effects are significant after FDR correction: the KS test gives $D = 0.515$ with $p_{\mathrm{FDR}} = 6.4 \times 10^{-3}$, and the Levene test gives $F = 9.74$ with $p_{\mathrm{FDR}} = 2.5 \times 10^{-3}$ (Tables~\ref{tab:ks_test} and~\ref{tab:levene_test}). The small size of the unknotted comparator should nevertheless be kept in mind when interpreting this comparison.

For the figure-eight slipknot $S4(1)$ ($n=127$), compared against $48$ unknotted homologs, the KDE in Fig.~\ref{fig:figure_eight_families_kde}\subref{fig:S4(1)_Curv_kde_median} shows a narrower and more negative distribution for the slipknotted family. The $S4(1)$ values are concentrated between $-30$ and $-22$ with a peak near $-24$, while the unknotted homologs are more dispersed. The two groups share the same median of $-26$, but their spreads differ substantially: $\mathrm{SD}_{S4(1)} = 2.08$ versus $\mathrm{SD}_{\mathrm{unknot}} = 6.04$, giving a variance ratio of $0.12$. This variance difference is highly significant after false discovery rate (FDR) correction (Levene test: $F = 18.38$, $p_{\mathrm{FDR}} = 4.0 \times 10^{-5}$; Table~\ref{tab:levene_test}), whereas the Kolmogorov-Smirnov (KS) test does not detect a significant global difference between the full distributions ($D = 0.178$, $p_{\mathrm{FDR}} = 0.190$; Table~\ref{tab:ks_test}). Thus, for $S4(1)$, the main signal lies in the tighter concentration of curvature values rather than in a median shift.

Overall, the figure-eight families follow the same broad pattern as the trefoil families, but with a stronger contrast for $K4(1)$ than for $S4(1)$.

\subsection{Summary of Statistical Tests across all Four Families}

Table~\ref{tab:ks_test} and Table~\ref{tab:levene_test} in Appendix~\ref{appendix c} report the full numerical results. The key findings are summarised here.

\paragraph{Levene test (variance).}
All four knotted and slipknotted families exhibit significantly lower variance in median Forman-Ricci curvature than their unknotted homologs after FDR correction ($p_{\text{FDR}} \leq 4.0 \times 10^{-3}$ for all four). Levene $F$-statistics range from $9.74$ ($K4(1)$) to $75.53$ ($S{+}3(1)$), and variance ratios range from $0.120$ ($S4(1)$) to $0.432$ ($S{+}3(1)$), all consistently below unity. This confirms that topological entanglement whether full knotting or slipknotting is universally associated with a narrower, more constrained curvature regime across all four families.

\paragraph{Kolmogorov-Smirnov test (distributional shift).}
Three of the four families reach statistical significance after FDR correction: $K4(1)$ ($D = 0.515$, $p_{\text{FDR}} = 6.41 \times 10^{-3}$), $K{+}3(1)$ ($D = 0.167$, $p_{\text{FDR}} = 2.51 \times 10^{-2}$), and $S{+}3(1)$ ($D = 0.214$, $p_{\text{FDR}} = 1.59 \times 10^{-10}$). The exception is $S4(1)$ ($D = 0.178$, $p_{\text{FDR}} = 1.90 \times 10^{-1}$), which does not show a significant global distributional difference despite exhibiting the most extreme variance contrast of the four families. This dissociation is not contradictory: the KS test is sensitive to differences in the full cumulative distribution function (CDF) (encompassing location, shape, and spread), whereas the Levene test isolates variance only. For $S4(1)$, the two groups share the same median, so the CDF-wide comparison is not significant, but the variance difference is and hence the curvature constraint imposed by slipknotting is clearly present and significant.

\subsection{Validation on Computationally Generated Random Loops}
\label{subsec:loops}

To test whether the curvature signal extends beyond the protein setting, we applied the pipeline to synthetic random polymer loops generated with \texttt{Topoly}~\cite{dabrowski2021topoly} at chain lengths $L \in \{100, 150, \cdots,  500\}$, with 500 loops per topology class per length. This provides a controlled setting in which sequence and fold-family effects are absent, and allows us to ask whether the observed signal is a more general consequence of topological entanglement rather than something specific to proteins.

\begin{figure}[ht]
    \centering
    \includegraphics[width=0.8\linewidth,height=7cm]{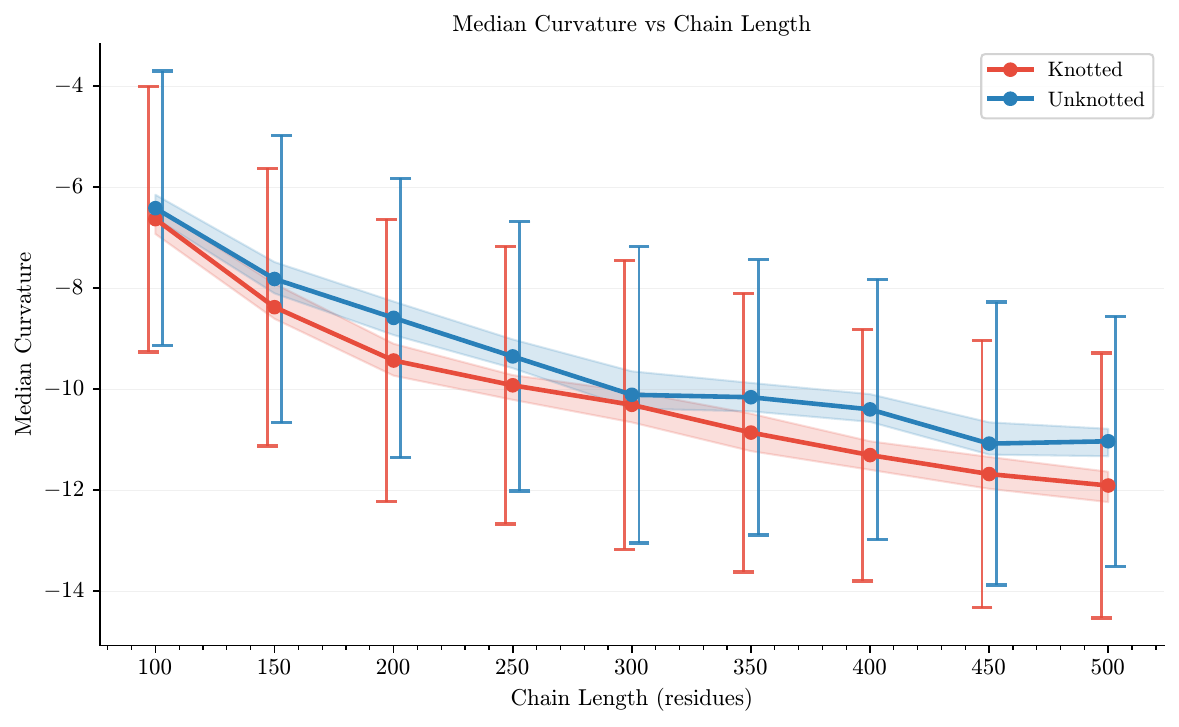}
    \caption{Per-loop median undirected Forman-Ricci curvature as a function of chain length for random knotted (red) and unknotted (blue) loops generated using \texttt{Topoly}. For each chain length, the plotted value is the median of the 500 per-loop medians; shaded regions indicate the 95\% bootstrap confidence interval. The two groups are nearly indistinguishable at $L = 100$. A detectable separation emerges at $L = 150$ and grows with chain length. Across lengths, the knotted group consistently exhibits substantially more negative curvature than the unknotted group.}
    \label{fig:line_median}
\end{figure}

\paragraph{Chain-length dependence of the curvature separation.}
Across all chain lengths, knotted loops display more negative median Forman-Ricci curvature than unknotted loops. The separation is small at $L = 100$ and becomes better detectable at $L = 150$. The separation grows substantially from $L = 200$ onward, with a transient narrowing at $L = 300$ where the two median lines nearly converge before diverging again. From $L = 400$ onward the separation is large and stable (Fig.~\ref{fig:line_median}), indicating that the curvature signal is a robust feature of knotted topology at longer chain lengths.

\paragraph{Statistical significance of the separation.}
The KS test was applied at each chain length with FDR correction applied simultaneously across all nine lengths. The test is not significant at $L = 100$ or $L = 150$ ($p_{\text{FDR}} > 0.05$; Fig.~\ref{fig:stat}), consistent with the visual impression from the trend plot. At all remaining lengths, from $L = 200$ onward, the distributional difference is statistically significant, confirming that the curvature signal becomes reliably detectable as chain length increases and topological constraints accumulate.

\subsection{Summary}

Across all four protein families and the synthetic loop validation, the median undirected Forman-Ricci curvature of knotted and slipknotted structures is consistently more negative than that of their unknotted counterparts. All four protein families exhibit a significant reduction in curvature variance relative to their comparator groups (Levene test, FDR-corrected), and three of the four show a significant global distributional difference (KS test, FDR-corrected). The exception, $S4(1)$, shows a non-significant KS result but the most extreme variance contrast of all four families; this is consistent with the hypothesis that slipknotting constrains the curvature regime without necessarily displacing the distribution's centre. The synthetic loop analysis confirms that these patterns are recoverable from topological structure alone, independently of protein biology, and that the group-level curvature separation strengthens as chain length increases. This is consistent with the increase in the occurrence probability of complex knots as the length increases. 

Taken together, these observations support the interpretation that knotting increases the overlap between persistent cycles in the hypergraph, raising the constituent hyperedge degrees and therefore lowering $F(e) = 2|e| - D$. This tighter structural organisation is consistent with the increased topological constraints imposed by chain entanglement, and it is recovered by the proposed pipeline across different knot types, slipknot configurations, and chain lengths.

\section{Conclusions}

The method developed here suggests that knotting leaves a detectable signature in persistent-homology-derived curvature statistics. Across the four protein families analysed, $K4(1)$, $S4(1)$, $K{+}3(1)$, and $S{+}3(1)$, knotted and slipknotted proteins consistently show more negative median Forman--Ricci curvature and lower variance than their unknotted homologs. Taken together with the synthetic-loop experiments, these results support the more general picture that topological entanglement is reflected in the organisation of persistent cycles, and can be recognised through simple curvature-based descriptors on the associated hypergraph.

There are several natural directions for further work. On the applied side, these descriptors could be incorporated into supervised learning pipelines as topologically informed features for recognising knot type or related forms of entanglement. More broadly, the same construction could be extended to other topological motifs, such as links or lassos, and to higher-dimensional persistence, where additional generators may capture further aspects of geometric confinement. More generally, the results suggest that persistent homology can provide simple and interpretable signatures of knotting beyond the specific protein setting that motivated this study.

  \paragraph{\textcolor{purple!60}{Authors' contributions.}}
  A.J.K.: conceptualisation, software (pipeline implementation), formal analysis, writing original draft, writing \& editing, validation.
  A.B.: conceptualisation, software (computational contributions), writing original draft, writing \& editing, validation, supervision.
  C.F.: conceptualisation, writing \& editing, validation, supervision.
  C.A.A.: reviewing, validation, supervision.
  All authors read and approved the final manuscript.

\paragraph{\textcolor{purple!60}{Data accessibility.}}
  Protein structure coordinates were obtained from the Protein Data Bank (\url{https://www.rcsb.org}).
  Knot annotations and homolog assignments were retrieved from KnotProt~2.0 (\url{https://knotprot.cent.uw.edu.pl}).
  Synthetic polymer loops were generated using the \texttt{Topoly} Python package.
  All code for the pipeline and statistical analyses is available at
  \url{https://github.com/aureliejodelle/PHyperRicci}.

  \paragraph{\textcolor{purple!60}{Funding.}}
  The authors acknowledge support from the Carnegie Corporation of New York through the AIMS Research and Innovation Centre.

  \paragraph{\textcolor{purple!60}{Acknowledgements.}}
  This publication was made possible by a grant from Carnegie Corporation of New York (provided through the AIMS Research and Innovation Centre).
  The statements made and views expressed are solely the responsibility of the author(s).

 \appendix

 \newpage
              
  \section*{Appendix A: Data Preparation}\label{appendix a}     
  Protein IDs and chain identifiers were retrieved programmatically from KnotProt. The analysis is restricted to four classes ($K{+}3(1)$, $S{+}3(1)$, $K4(1)$, $S4(1)$), as they are the most extensively catalogued knot types in the database.                            
  The comparative analysis retains the KnotProt-annotated knotted proteins and their unknotted sequence-similar homologs, since the objective is to contrast topological features of knotted proteins with those of structurally related but unknotted counterparts.              
                                                     
  \paragraph{Classification of unannotated PDB homologs.}   
  Proteins in the PDB homolog set lack knot-type annotations in KnotProt and therefore require independent classification prior to inclusion. We classify them into knotted or unknotted by computing knot polynomials after taking stocastic closures using Python's software \texttt{Topoly}. Proteins classified as knotted were excluded; the remaining unknotted proteins were merged with the directly retrieved unknotted homologs to form a unified unknotted homolog set for each class. 
  
  \paragraph{Homolog group sizes.}                          
  Unknotted families are systematically smaller than the corresponding knotted classes.                                                              
  This asymmetry arises because a single unknotted protein may be homologous to several knotted proteins simultaneously, so the number of unique unknotted homologs is considerably lower than the number of knotted-unknotted pairs.                                                  
  For example, the human carbonic anhydrase chain $4\mathrm{WW8}_A$ is sequence-similar to nine distinct knotted proteins across the $K{+}3(1)$ class~\cite{dabrowski2019knotprot,jamroz2015knotprot}. 
 
  \paragraph{Synthetic random loop generation.}             
  Loops were generated with \textit{generate\_loop(L, 1, output='list')}, which produces a random walk of $L$ nodes with unit step size. Each knot type was determined by the Jones polynomial. For each length, ten independent sampling rounds were conducted. Within each round, loops were generated one at a time and accepted to either the knotted or unknotted pool until each pool had accumulated exactly 50 loops, yielding 500 loops per topology class per chain length and ensuring strictly balanced classes by construction. Loops and their PH data were stored in JSON format for downstream analysis.               
 
  \section*{Appendix B: Computational Pipeline}\label{appendix b} 
  Persistent homology of dimension one was computed using Ripserer on a point cloud with the exact call:                                         
  \begin{verbatim}                       
    PH = ripserer(grid; dim_max=1, alg=:involuted)     
  \end{verbatim}                        
  where \textit{grid} is a vector of three-element tuples $(x_i, y_i, z_i)$ constructed from the $C_\alpha$ coordinate CSV file.                                                   
  The flag \textit{alg=:involuted} activates the involuted algorithm for cycle representative extraction; \textit{dim\_max=1} restricts computation to $H_0$ and $H_1$.                       
  The computation runs in parallel across all proteins of a class using Julia's \textit{@threads} macro.                                    
  For each protein, three JSON outputs are written: \textit{PH\_1/(id).json} (full barcodes and representatives), \textit{representatives/(id).json} (cycles only, used as hyperedge input), and \textit{barcodes/(id).json} (birth-death pairs only, used for visualisation). Infinite deaths are serialised as JSON \textit{null}.  
  \paragraph{Rationale for the Vietoris-Rips filtration.}                       
  The VR construction is the canonical choice for two reasons.                  
  First, it matches the filtration used in hyperTDA~\cite{barbensi2022hypergraphs}, enabling a direct comparison of our hyperedge construction with the original framework. Second, and more critically, we do not impose any threshold on the filtration: all $H_1$ generators produced up to the point cloud diameter are retained. This design reflects our aim of capturing the \emph{global} topological organisation of the protein backbone. Even short-lived cycles contribute to the hyperedge degree count $D$ and therefore influence the Forman-Ricci curvature of neighbouring hyperedges. The alpha complex, whilst computationally more efficient and geometrically better adapted to three-dimensional data, enforces geometric admissibility constraints absent from the VR filtration. Since our downstream analysis depends on the completeness of the cycle set rather than on distances between persistence diagrams, the VR filtration is the appropriate construction here. A practical consequence is that for large proteins, lots of short-lived generators are retained; these contribute to hyperedge degree sums and influence $F(e)$, but their individual geometric meaning is limited.                                              
  \paragraph{Cycle representatives.}  
  For each $H_1$ generator returned by \texttt{Ripserer.jl}, the involuted algorithm provides a representative cycle stored as a
  collection of ordered pairs of $C_\alpha$ atom indices $[v_1, v_2]$. Concretely, for a generator $c$:                                              
  \begin{equation}                    
      \mathrm{rep}(c) =                         
      \bigl\{[v_1, v_2] \;\colon\; [v_1, v_2] \in c \bigr\}.                    
  \end{equation}                                                              
  Representatives are used as returned by \texttt{Ripserer.jl} without further post-processing.
  The involuted algorithm satisfies a boundary-matrix symmetry condition that biases representatives towards minimal cycles whilst remaining computationally tractable at the scales required for proteins with hundreds to thousands of residues~\cite{cufar2020ripserer,bauer2021ripser}; it therefore produces structurally more interpretable cycles than generic pivot-based algorithms, without incurring the cost of explicit length-minimisation.

  \paragraph{Hypergraph construction.}
  The hypergraph is constructed from cycle representatives. Node indices equal to zero or negative are discarded as artefacts of the \texttt{Ripserer.jl} indexing convention; all valid indices are one-based residue positions. The hyperedge map for each protein is serialised as \textit{hyperedge\_map/(id).json}, a dictionary mapping integer hyperedge IDs to sorted lists of node indices. Processing is parallelised over proteins within each class using eight concurrent threads.

  \section*{Appendix C: Statistical Analysis}\label{appendix c}

  To assess whether the median Forman-Ricci curvature distinguishes knotted proteins from their unknotted homologs, we employ two complementary non-parametric tests applied to each of the four class pairs. No distributional assumptions are made about the curvature values.

  \paragraph{Kolmogorov-Smirnov test.}
  The two-sample KS test assesses whether the two groups are drawn from the same underlying distribution. The null hypothesis $H_0$ is that the empirical distributions of median curvature in the knotted class and the unknotted homolog set are identical. The KS $D$-statistic measures the maximum absolute difference between the two empirical cumulative distribution functions and ranges from 0 to 1; larger values indicate greater distributional separation.
  Group pairs with fewer than three observations are excluded from this test.

  \paragraph{Levene test.}
  The Levene test (computed via \texttt{scipy.stats.levene},
  \texttt{center='median'}) assesses equality of variance between the two groups.
  Using the median as the centring statistic makes the test robust to non-normality and is the Brown-Forsythe variant, which is the standard choice for non-normal or skewed distributions.
  The null hypothesis $H_0$ is that the two groups have equal population variances.
  This test is included because variance differences, in addition to location shifts, are biologically meaningful: a knotted class with substantially lower variance than its unknotted counterpart suggests that topological entanglement imposes a structural constraint on the curvature distribution, irrespective of any shift in its centre. The variance ratio $\mathrm{Var}(\text{knotted})/\mathrm{Var}(\text{unknotted})$ is reported alongside the Levene $F$-statistic and $p$-value as an interpretable effect size.

  \paragraph{Multiple testing correction.}
  All raw $p$-values are corrected for multiple comparisons using the Benjamini-Hochberg false discovery rate procedure, implemented via the command \\ \texttt{statsmodels.stats.multitest.multipletests} with \textit{alpha=0.05} and \texttt{method='fdr\_bh'}. The FDR correction is applied \emph{separately} for the KS and Levene tests, each correction spanning all four class comparisons simultaneously. Results are reported at the $\alpha = 0.05$ level after correction. For the synthetic loop validation, the KS correction spans all nine chain lengths simultaneously.

  \paragraph{Chain length as a potential confound.}
  The number of $C_\alpha$ atoms varies across proteins and could in principle influence the number of persistent generators, hyperedge degree sums, and therefore median curvature, independently of topology.
  Chain-length distributions for all eight groups are reported in Table~\ref{tab:dataset_methods_compact2}.
  In three out of four families, the unknotted homologs have median chain lengths that are comparable to or longer than the knotted proteins: $S4(1)$ ($507$ vs $537$), $K4(1)$ ($478$ vs $462$),
  $K{+}3(1)$ ($257$ vs $310$), and $S{+}3(1)$ ($259$ vs $427$).
  If longer chains produced more negative curvature, the unknotted groups would exhibit the same or more negative median curvature than the knotted groups.
  The observed direction is the opposite in every family, indicating that the chain-length effect acts against the observed curvature signal rather than in favour of it, providing indirect evidence that the curvature difference is topological rather than length-driven.

  \begin{table}[ht]
  \centering
  \footnotesize
  \caption{\textbf{Table~A1.}
  Levene test of variance equality between each knotted/slipknotted protein family and its unknotted homolog comparator group, applied to the per-protein median Forman-Ricci curvature (\textit{Curv\_median}).
  The test statistic is the standard Levene $F$-statistic centred on group medians, robust to non-normality. The null hypothesis is that the two groups have equal population variance. The variance ratio $\mathrm{Var}(\text{knotted})/\mathrm{Var}(\text{unknotted})$ quantifies the relative dispersion; values below 1 indicate the knotted group is more tightly concentrated. Raw and FDR-corrected $p$-values are both reported; the Benjamini-Hochberg correction is applied simultaneously across all four comparisons. All four families are significant after FDR correction at $\alpha=0.05$. Sample sizes: $K4(1)$ $n=70$ vs $13$; $S4(1)$ $n=127$ vs $48$; $K{+}3(1)$ $n=615$ vs $94$; $S{+}3(1)$ $n=1{,}211$ vs $340$. $^{\dagger}$For $K4(1)$, the raw and FDR-corrected $p$-values are identical because this is the largest (least significant) of the four $p$-values; the Benjamini-Hochberg multiplier for rank $i=4$ out of $n=4$ comparisons is $n/i=1$, so the adjusted value equals the raw value.}
  \label{tab:levene_test}
  \begin{tabular}{l r r r r l}
  \toprule
  \textbf{Comparison} &
  \textbf{Levene $F$} &
  \textbf{$p$ (raw)} &
  \textbf{$p$ (FDR)} &
  \textbf{Var.\ ratio} &
  \textbf{Sig.} \\
  \midrule
  $K4(1)$ vs unknotted
      & 9.74  & $2.50\times10^{-3}$ & $2.50\times10^{-3}$$^{\dagger}$
      & 0.189 & \cellcolor{green!20}Yes~(**) \\
  $S4(1)$ vs unknotted
      & 18.38 & $3.00\times10^{-5}$ & $4.00\times10^{-5}$
      & 0.120 & \cellcolor{green!20}Yes~(****) \\
  $K{+}3(1)$ vs unknotted
      & 51.19 & $2.10\times10^{-12}$ & $4.20\times10^{-12}$
      & 0.350 & \cellcolor{green!20}Yes~(****) \\
  $S{+}3(1)$ vs unknotted
      & 75.53 & $8.98\times10^{-18}$ & $3.59\times10^{-17}$
      & 0.432 & \cellcolor{green!20}Yes~(****) \\
  \bottomrule
  \end{tabular}
  \end{table}

  \begin{table}[ht]
  \centering
  \footnotesize
  \caption{\textbf{Table~A2.}
  Two-sample Kolmogorov-Smirnov (KS) test comparing the per-protein median Forman-Ricci curvature distributions of each knotted/slipknotted family against its unknotted homolog comparator. The KS $D$-statistic measures the maximum absolute difference between the two empirical CDFs.
  Raw and FDR-corrected $p$-values are both reported; the Benjamini-Hochberg correction is applied simultaneously across all four comparisons and separately from the Levene corrections in Table~A1. Significance: \textit{ns} $p\geq0.05$; * $p<0.05$; ** $p<0.01$; **** $p<0.0001$ (all FDR-corrected). The $S4(1)$ family does not reach significance, attributable to the relatively small unknotted comparator ($n=48$).}
  \label{tab:ks_test}
  \begin{tabular}{l r r r r l}
  \toprule
  \textbf{Comparison} &
  \textbf{KS $D$} &
  \textbf{Med.\ diff.} &
  \textbf{$p$ (raw)} &
  \textbf{$p$ (FDR)} &
  \textbf{Sig.} \\
  \midrule
  $K4(1)$ vs unknotted
      & 0.515 & $-3$ & $3.20\times10^{-3}$ & $6.41\times10^{-3}$
      & \cellcolor{green!20}Yes~(**) \\
  $S4(1)$ vs unknotted
      & 0.178 & $0$  & $1.90\times10^{-1}$ & $1.90\times10^{-1}$
      & \cellcolor{red!20}No~(ns) \\
  $K{+}3(1)$ vs unknotted
      & 0.167 & $-1$ & $1.88\times10^{-2}$ & $2.51\times10^{-2}$
      & \cellcolor{green!20}Yes~(*) \\
  $S{+}3(1)$ vs unknotted
      & 0.214 & $-1$ & $3.96\times10^{-11}$ & $1.59\times10^{-10}$
      & \cellcolor{green!20}Yes~(****) \\
  \bottomrule
  \end{tabular}
  \end{table}

  \begin{figure}[ht]
      \centering
      \includegraphics[width=0.9\linewidth]{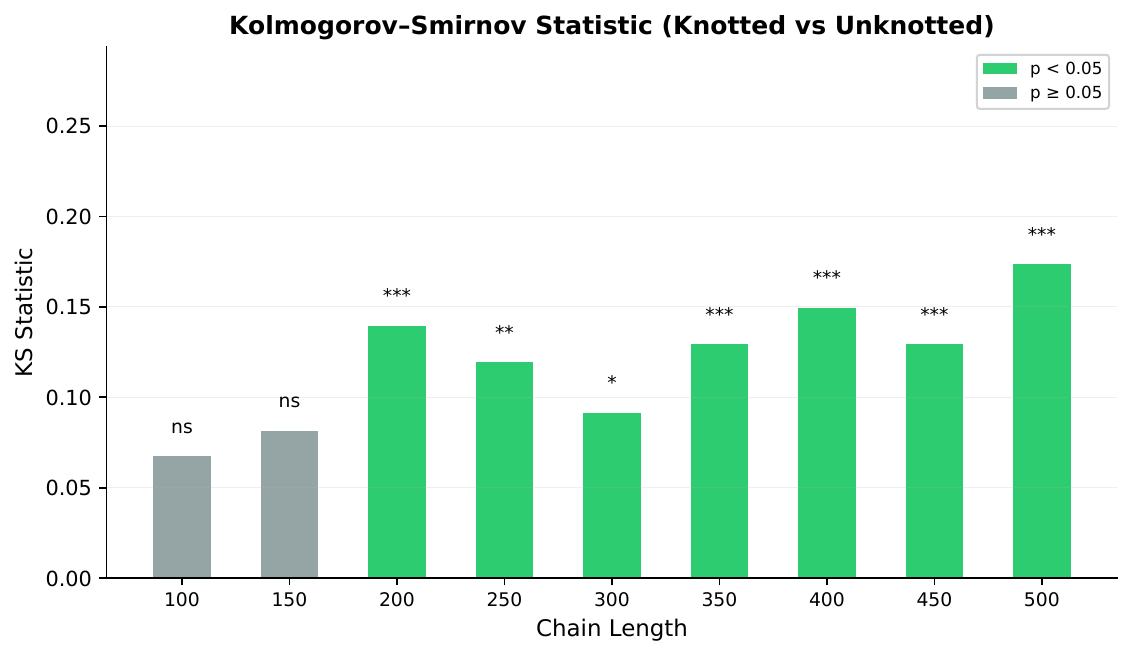}
      \caption{Kolmogorov-Smirnov $D$-statistic for the comparison of knotted versus unknotted random loops at each chain length ($L=100,150,\ldots,500$). Bars are coloured green where the FDR-corrected $p$-value is $<0.05$ and grey otherwise. Significance: \textit{ns} $p\geq0.05$, * $p<0.05$, ** $p<0.01$, *** $p<0.001$ (Benjamini-Hochberg, applied simultaneously across all nine chain lengths).
      The separation is not significant at $L=100$ and $L=150$, consistent with the trend plot in Fig.~\ref{fig:line_median}.}
      \label{fig:stat}
  \end{figure}

  \end{document}